\documentclass{amsart}
\usepackage[margin=1in]{geometry} 
\usepackage[backend=biber, style=alphabetic]{biblatex}
\usepackage{mymacros, pgf}
\usetikzlibrary{calc,arrows}
\usepackage{subcaption}
\usepackage[foot]{amsaddr}
\bibliography{ternbib} % or
\title{Ternary arithmetic, factorization, and the class number one problem}
\author{Aram Bingham}
\address{Tulane University, 6823 St. Charles Ave., New Orleans, LA, 70118}
\email{abingham@tulane.edu}
\date{Feb. 2020}
\subjclass[2010]{11A05, 11H06, 11Y05}

\numberwithin{equation}{section}

\begin{document}

\begin{abstract}
Ordinary multiplication of natural numbers can be generalized to a ternary operation by considering discrete volumes of lattice hexagons. With this operation, a natural notion of `3-primality' -- primality with respect to ternary multiplication -- is defined, and it turns out that there are very few 3-primes. They correspond to imaginary quadratic fields $\Q(\sqrt{-n})$, $n>0$, with odd discriminant and whose ring of integers admits unique factorization. We also describe how to determine representations of numbers as ternary products and related algorithms for usual primality testing and integer factorization.


\end{abstract}

\maketitle\vspace{-3em}
\section{Basic Ideas}
When ideas become engrained, it can be hard to imagine other possibilities. But in escaping from deep-seated notions, we may be surprised by what we find. In this spirit, this article will present a deformation of integer arithmetic which remains grounded in geometry and leads to new perspectives on old problems in number theory related to primality and factorization. 

Let's start with the absolute basics. Say you want to multiply two whole numbers, $a$ and $b$. What do you do to find the product $ab$? 

One option is the following. Draw $a$ parallel lines on a piece of paper. Now draw $b$ lines which are parallel to each other but perpendicular to the $a$ parallel lines you drew first. The number of intersection points of the two sets of lines is your product $ab$. Further, one sees the commutativity of multiplication in the fact that the number of intersection points doesn't change when you rotate the whole picture by $90^\circ$.
\begin{center}
\begin{figure}[htbp!]
\begin{tikzpicture}[scale=0.5]
\foreach \x in {1, 2, 3, 4, 5, 6,7}
	\draw (\x,-0.5) --(\x, 4.5) ;
\foreach \y in {1,2,3}
	\draw (-0.5, \y) -- (8.5, \y);
\foreach \x in {1, 2, 3, 4, 5, 6,7}	
	\foreach \y in {1,2,3}	
	\filldraw[color=red, fill=red] (\x, \y) circle (1mm);
\end{tikzpicture} \qquad
\begin{tikzpicture}[scale=0.5]
\filldraw[green!10!white] (0,0) rectangle (6,2);
\draw[gray, very thin] (-1.5,-1.5) grid (8.5,3.5);
\draw[<->] (-1.5,0) -- (8.5, 0);
\draw[<->] (0,-1.5) -- (0, 3.5);
%\foreach \x in {1, 2, 3, 4, 5, 6,7}
	%\draw (\x,-1) --(\x, 5) ;
%\foreach \y in {1,2,3}
	%\draw (-1, \y) -- (9, \y);
\foreach \x in {0,1, 2, 3, 4, 5, 6}	
	\foreach \y in {0,1,2}	
	\filldraw[color=red, fill=red] (\x, \y) circle (1mm);
\end{tikzpicture}
\caption{3 times 7 is 21}
\label{fig1}
\end{figure}
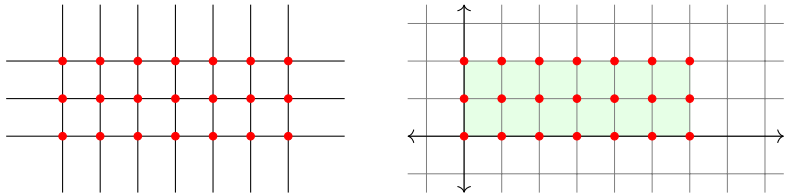
\end{center}
\vspace{-3em}
Let's think of this slightly more formally. Given a lattice of points in the plane, we will define a \textbf{lattice polygon} to be a polygon whose vertices are lattice points and whose edges are only in the directions of nearest neighbors from a given vertex. Using the $\Z^2$ lattice, this means that edges are either in the horizontal or vertical directions. The product $ab$ is then realized as the number of lattice points inside or on the boundary of the lattice rectangle with corners at $(0,0)$, $(b-1,0)$, $(0,a-1)$ and $(a-1, b-1)$. Restated, $ab$ is the \textbf{discrete volume} of the lattice rectangle with $a$ points along two opposite edges in the vertical direction and $b$ lattice points along the opposite edges in the horizontal direction; see the right side of Figure~\ref{fig1}. This view of multiplication allows us to codify the following simple observation.

\begin{nfact} \label{fact1} A number is prime if and only if it cannot be represented as the discrete volume of a $\Z^2$ lattice rectangle (with edges in the horizontal and vertical directions) and where each edge contains at least two lattice points.
\end{nfact}

In this model, the commutativity of multiplication is beheld in the preservation of discrete volume when interchanging which lattice direction corresponds to which factor in the product $ab$. This suggests that sensible alternatives to standard multiplication might then be found in by taking discrete volumes of other lattice polygons in other lattices.  

Opting for maximal symmetry, we consider the \textbf{hexagonal lattice} \cite[pp.~60-61]{coxeterGeometry}. In this lattice, rotation of the plane by $60^\circ$ about any lattice point sends lattice points to other lattice points, arranging the nearest neighbors of any lattice point $P$ in a regular hexagon. These six other points come in pairs along three lines through $P$ as compared to the four nearest neighbors of a point in the $\Z^2$ lattice which come in pairs along two lines (Figure~\ref{fig:neighbors}). Hence we gain an extra direction that can be assigned to a third factor.

\begin{center}
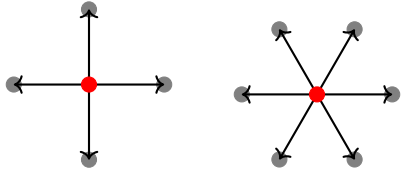
\begin{figure}[htbp!]
\begin{tikzpicture}
	\filldraw[color=gray] (1,0) circle (1mm) (-1, 0) circle (1mm) (0, 1) circle (1mm) (0, -1) circle (1mm);
	\draw [<->, thick] (0,-1) --(0, 1) ;
	\draw [<->, thick] (-1,0) --(1, 0) ;
	\filldraw[color=red, fill=red] (0, 0) circle (1mm);
	\end{tikzpicture} \qquad
\begin{tikzpicture}
\coordinate (0;0) at (0,0); 
\foreach \c in {1,...,3}{%  
\foreach \i in {0,...,5}{% 
\pgfmathtruncatemacro\j{\c*\i}
\coordinate (\c;\j) at (60*\i:\c);  
} }
\foreach \i in {0,...,5}{% 
\pgfmathtruncatemacro\j{\i}
\filldraw[color=gray] (1;\j) circle (1mm);}
	\draw [<->, thick] (1;0) --(1;3) ;
	\draw [<->, thick] (1;1) --(1;4) ;
	\draw [<->, thick] (1;2) --(1;5) ;
\filldraw[color=red, fill=red] (0, 0) circle (1mm);	
\end{tikzpicture}
\caption{Nearest neighbors in a square ($\Z^2$) lattice and a hexagonal lattice. }
\label{fig:neighbors}
\end{figure}
\end{center}
\vspace{-3em}

Recall that the \textbf{arity} of a function or algebraic operation refers to the number of inputs or arguments it takes. Binary operations can always be iterated to fabricate operations of higher arity, but we will introduce a true \textbf{ternary product} on the natural numbers.\footnote{Throughout this manuscript, we take the set of natural numbers $\N$ to start at 1.} While distinct from repeated multiplication, it bears some of the same nice properties: commutativity guaranteed by geometry, and the presence of a \emph{multiplicative identity} in the number 1. We will denote this product as
\[\la -, -,-\ra: \N^3 \to \N,\]
and define it, in analogy with our lattice model of binary multiplication, as the function which takes the triplet $(a,b,c)$ to the number of lattice points inside the equiangular lattice hexagon with $a$ points along two opposite edges, $b$ points along the next pair of edges, and $c$ points along the final pair. Illustrations are given in Figures~\ref{fig:ternsym} and~\ref{fig:ternbin}. 
\begin{center}
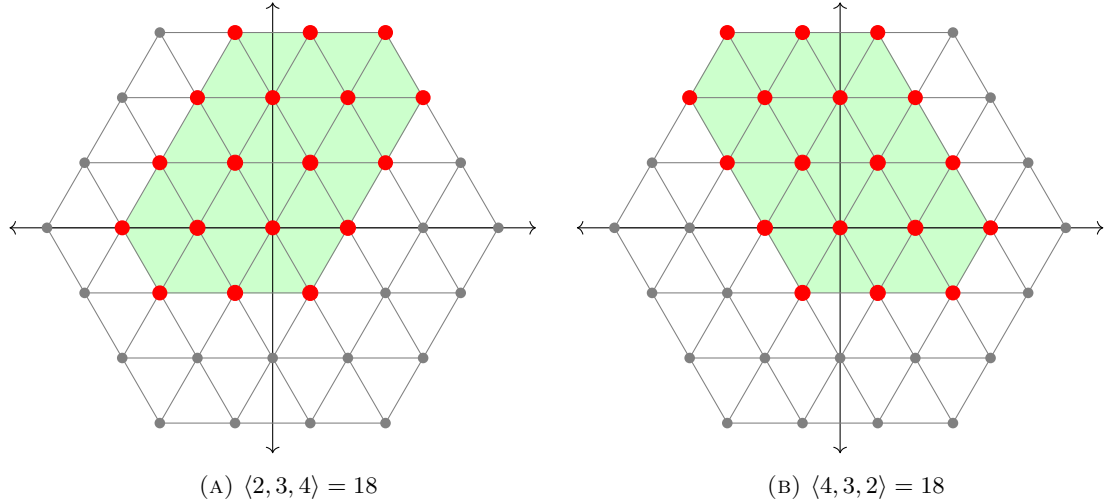
\begin{figure}[htbp!]
\begin{subfigure}[b]{0.45\textwidth}
\begin{tikzpicture}
%%%  define vertices with coordinates
\coordinate (0;0) at (0,0); 
\foreach \c in {1,...,3}{%  
\foreach \i in {0,...,5}{% 
\pgfmathtruncatemacro\j{\c*\i}
\coordinate (\c;\j) at (60*\i:\c);  
} }
\foreach \i in {0,2,...,10}{% 
\pgfmathtruncatemacro\j{mod(\i+2,12)} 
\pgfmathtruncatemacro\k{\i+1}
\coordinate (2;\k) at ($(2;\i)!.5!(2;\j)$) ;}

\foreach \i in {0,3,...,15}{% 
\pgfmathtruncatemacro\j{mod(\i+3,18)} 
\pgfmathtruncatemacro\k{\i+1} 
\pgfmathtruncatemacro\l{\i+2}
\coordinate (3;\k) at ($(3;\i)!1/3!(3;\j)$)  ;
\coordinate (3;\l) at ($(3;\i)!2/3!(3;\j)$)  ;
 }

%%shape
\filldraw[green!20!white] (2;6) -- (2;7) -- (1;5) -- (3;2) -- (3;3) -- (3;5) ;
 \draw[<->] (0,-3) -- (0, 3);
 
 %%%%%%%%% draw lines %%%%%%%%
 \foreach \i in {0,...,6}{% 
 \pgfmathtruncatemacro\k{\i}
 \pgfmathtruncatemacro\l{15-\i}
 \draw[thin,gray] (3;\k)--(3;\l);
 \pgfmathtruncatemacro\k{9-\i} 
 \pgfmathtruncatemacro\l{mod(12+\i,18)}   
 \draw[thin,gray] (3;\k)--(3;\l); 
 \pgfmathtruncatemacro\k{12-\i} 
 \pgfmathtruncatemacro\l{mod(15+\i,18)}   
 \draw[thin,gray] (3;\k)--(3;\l);}    
 \draw[<->] (-3.5,0) -- (3.5, 0);
%%% draw gray circles at lattice points
\fill [gray] (0;0) circle (2pt);
 \foreach \c in {1,...,3}{%
 \pgfmathtruncatemacro\k{\c*6-1}    
 \foreach \i in {0,...,\k}{% 
   \fill [gray] (\c;\i) circle (2pt);}}  
%%%%% points in hexagon
\foreach \x in {0,1, 2, 3, 4, 5}{%
	\pgfmathtruncatemacro\y{\x}	
	\filldraw[color=red, fill=red] (1;\y) circle (1mm);}
\foreach \x in {1,...,7}	{ \pgfmathtruncatemacro\y{\x}
	\fill [red] (2;\x) circle (1mm);}	
\foreach \x in {2,...,5}	{ \pgfmathtruncatemacro\y{\x}
	\fill [red] (3;\x) circle (1mm);}		
\fill [red] (0;0) circle (1mm);
\end{tikzpicture} 
 
\caption{$\la 2, 3, 4\ra = 18$ }
\end{subfigure}
\begin{subfigure}[b]{0.45\textwidth}
\begin{tikzpicture}
%%%  define vertices with coordinates
\coordinate (0;0) at (0,0); 
\foreach \c in {1,...,3}{%  
\foreach \i in {0,...,5}{% 
\pgfmathtruncatemacro\j{\c*\i}
\coordinate (\c;\j) at (60*\i:\c);  
} }
\foreach \i in {0,2,...,10}{% 
\pgfmathtruncatemacro\j{mod(\i+2,12)} 
\pgfmathtruncatemacro\k{\i+1}
\coordinate (2;\k) at ($(2;\i)!.5!(2;\j)$) ;}
\foreach \i in {0,3,...,15}{% 
\pgfmathtruncatemacro\j{mod(\i+3,18)} 
\pgfmathtruncatemacro\k{\i+1} 
\pgfmathtruncatemacro\l{\i+2}
\coordinate (3;\k) at ($(3;\i)!1/3!(3;\j)$)  ;
\coordinate (3;\l) at ($(3;\i)!2/3!(3;\j)$)  ;
 }

%%shape
\filldraw[green!20!white]  (2;11) -- (2;0) -- (3;4) -- (3;6) -- (3;7) --(1;4) ;
 \draw[<->] (0,-3) -- (0, 3);

 %%%%%%%%% draw lines %%%%%%%%
 \foreach \i in {0,...,6}{% 
 \pgfmathtruncatemacro\k{\i}
 \pgfmathtruncatemacro\l{15-\i}
 \draw[thin,gray] (3;\k)--(3;\l);
 \pgfmathtruncatemacro\k{9-\i} 
 \pgfmathtruncatemacro\l{mod(12+\i,18)}   
 \draw[thin,gray] (3;\k)--(3;\l); 
 \pgfmathtruncatemacro\k{12-\i} 
 \pgfmathtruncatemacro\l{mod(15+\i,18)}   
 \draw[thin,gray] (3;\k)--(3;\l);}    
 
\draw[<->] (-3.5,0) -- (3.5, 0);
%%%%% grey circles at lattice points
\fill [gray] (0;0) circle (2pt);
 \foreach \c in {1,...,3}{%
 \pgfmathtruncatemacro\k{\c*6-1}    
 \foreach \i in {0,...,\k}{% 
   \fill [gray] (\c;\i) circle (2pt);}}  
%%%%% red points in parallelogram
\foreach \x in {0,1, 2, 3, 4, 5}{%
	\pgfmathtruncatemacro\y{\x}	
	\filldraw[color=red, fill=red] (1;\y) circle (1mm);}
\foreach \x in {0,...,5}	{ \pgfmathtruncatemacro\y{\x}
	\fill [red] (2;\x) circle (1mm);}	
\foreach \x in {4,...,7}	{ \pgfmathtruncatemacro\y{\x}
	\fill [red] (3;\x) circle (1mm);}		
\fill [red] (0;0) circle (1mm);
\fill [red] (2;11) circle (1mm);
 
\end{tikzpicture}  
\caption{$\la 4,3,2 \ra=18$}
\end{subfigure}
\caption{Commutativity of ternary multiplication under reflection.}
\label{fig:ternsym}
\end{figure}
\end{center}
\vspace{-3em}

\section{Properties of Ternary Arithmetic}
A lattice hexagon representing the product $\la a, b, c\ra$ can be acted upon by any of the symmetries of the lattice. Under this action, any lattice direction can be taken to any other while the discrete volume of a lattice hexagon is always preserved, implying that $\la -, -, - \ra$ is fully commutative as a ternary product.

Notice that if we put a 1 into one of the arguments of this product, one of the pairs of sides of the hexagon degenerates to a single point and we instead have a parallelogram (see Figure~\ref{fig:ternbin}). The discrete volume of this parallelogram is then just the value of the binary product of the other two arguments, so we observe ordinary multiplication as a specialization of the ternary product. Further, if we have a 1 appearing twice as an argument, then $\la 1, 1, n\ra$ leads to just a row of $n$ points (with discrete volume $n$), showing that $1$ indeed behaves as a multiplicative identity.

By analogy with Fact~\ref{fact1}, we make the following definition.
\begin{ndefn} We will say that a natural number $n$ is \textbf{$3$-prime} if it can not be represented as the discrete volume of an equiangular lattice hexagon for which at least two pairs of opposing sides have at least two points. Equivalently, $n$ is $3$-prime if there is no choice of $x, y, z$ such that $\la x, y, z\ra=n$ other than $\la 1, 1, n\ra=n$ and permutations of these inputs.
\end{ndefn}\vspace{-1em}
To avoid confusion with the usual definition of primality, from now on we will say that a natural number $n$ is \textbf{2-prime} to mean that its only natural number factors are $1$ and $n$.\footnote{Under binary multiplication, of course.} We shall also say that a number is ``$2$-composite'' or ``$3$-composite'' to mean that it is not $2$-prime or not $3$-prime, respectively. An immediate consequence of this definition is that 3-primality implies 2-primality, but not vice versa. For instance \mbox{$\la 2, 2, 2\ra= 7$} is not $3$-prime, but 2, 3 and 5 are still $3$-prime, and with a little checking you can convince yourself that 11 is as well. This begs the following question.

\begin{center}
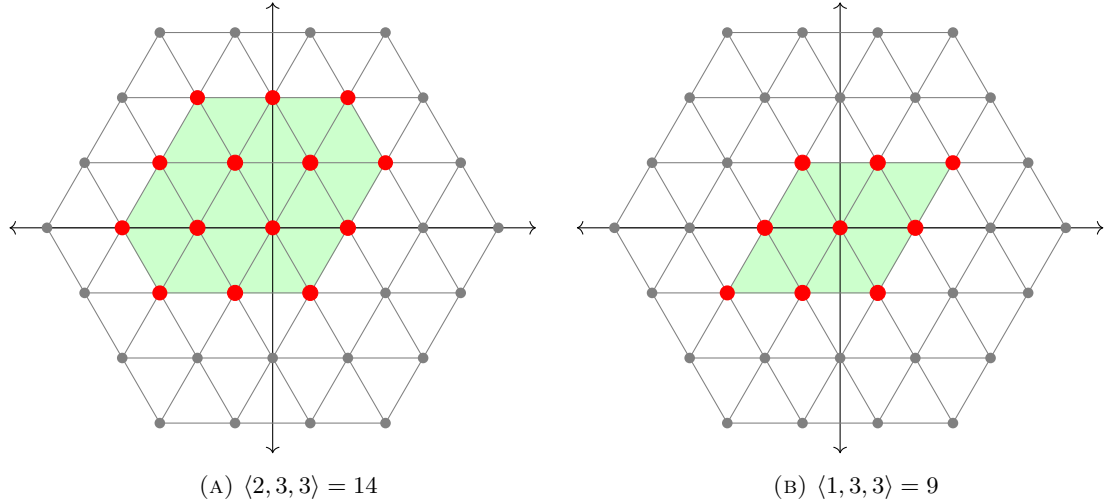
\begin{figure}[htbp!]
\begin{subfigure}[b]{0.45\textwidth}
\begin{tikzpicture}
%%%  define vertices with coordinates
\coordinate (0;0) at (0,0); 
\foreach \c in {1,...,3}{%  
\foreach \i in {0,...,5}{% 
\pgfmathtruncatemacro\j{\c*\i}
\coordinate (\c;\j) at (60*\i:\c);  
} }
\foreach \i in {0,2,...,10}{% 
\pgfmathtruncatemacro\j{mod(\i+2,12)} 
\pgfmathtruncatemacro\k{\i+1}
\coordinate (2;\k) at ($(2;\i)!.5!(2;\j)$) ;}

\foreach \i in {0,3,...,15}{% 
\pgfmathtruncatemacro\j{mod(\i+3,18)} 
\pgfmathtruncatemacro\k{\i+1} 
\pgfmathtruncatemacro\l{\i+2}
\coordinate (3;\k) at ($(3;\i)!1/3!(3;\j)$)  ;
\coordinate (3;\l) at ($(3;\i)!2/3!(3;\j)$)  ;
 }

%%shape
\filldraw[green!20!white] (2;6) -- (2;7) -- (1;5) -- (2;1) -- (2;2) -- (2;4) ;
 \draw[<->] (0,-3) -- (0, 3);
 
 %%%%%%%%% draw lines %%%%%%%%
 \foreach \i in {0,...,6}{% 
 \pgfmathtruncatemacro\k{\i}
 \pgfmathtruncatemacro\l{15-\i}
 \draw[thin,gray] (3;\k)--(3;\l);
 \pgfmathtruncatemacro\k{9-\i} 
 \pgfmathtruncatemacro\l{mod(12+\i,18)}   
 \draw[thin,gray] (3;\k)--(3;\l); 
 \pgfmathtruncatemacro\k{12-\i} 
 \pgfmathtruncatemacro\l{mod(15+\i,18)}   
 \draw[thin,gray] (3;\k)--(3;\l);}    
 
 \draw[<->] (-3.5,0) -- (3.5, 0);

%%% draw gray circles at lattice points
\fill [gray] (0;0) circle (2pt);
 \foreach \c in {1,...,3}{%
 \pgfmathtruncatemacro\k{\c*6-1}    
 \foreach \i in {0,...,\k}{% 
   \fill [gray] (\c;\i) circle (2pt);}}  
%%%%% points in hexagon
\foreach \x in {0,1, 2, 3, 4, 5}{%
\pgfmathtruncatemacro\y{\x}	
	\filldraw[color=red, fill=red] (1;\y) circle (1mm);}
\foreach \x in {1,...,7}	{ \pgfmathtruncatemacro\y{\x}
\fill [red] (2;\x) circle (1mm);}	
\fill [red] (0;0) circle (1mm);
\end{tikzpicture}  
\caption{$\la 2, 3, 3\ra = 14$ }
\end{subfigure}
\begin{subfigure}[b]{0.45\textwidth}
\begin{tikzpicture}
%%%  define vertices with coordinates
\coordinate (0;0) at (0,0); 
\foreach \c in {1,...,3}{%  
\foreach \i in {0,...,5}{% 
\pgfmathtruncatemacro\j{\c*\i}
\coordinate (\c;\j) at (60*\i:\c);  
} }
\foreach \i in {0,2,...,10}{% 
\pgfmathtruncatemacro\j{mod(\i+2,12)} 
\pgfmathtruncatemacro\k{\i+1}
\coordinate (2;\k) at ($(2;\i)!.5!(2;\j)$) ;}

\foreach \i in {0,3,...,15}{% 
\pgfmathtruncatemacro\j{mod(\i+3,18)} 
\pgfmathtruncatemacro\k{\i+1} 
\pgfmathtruncatemacro\l{\i+2}
\coordinate (3;\k) at ($(3;\i)!1/3!(3;\j)$)  ;
\coordinate (3;\l) at ($(3;\i)!2/3!(3;\j)$)  ;
 }

%%shape
\filldraw[green!20!white]  (2;7) -- (1;5) -- (2;1) -- (1;2) ;
 \draw[<->] (0,-3) -- (0, 3);

 %%%%%%%%% draw lines %%%%%%%%
 \foreach \i in {0,...,6}{% 
 \pgfmathtruncatemacro\k{\i}
 \pgfmathtruncatemacro\l{15-\i}
 \draw[thin,gray] (3;\k)--(3;\l);
 \pgfmathtruncatemacro\k{9-\i} 
 \pgfmathtruncatemacro\l{mod(12+\i,18)}   
 \draw[thin,gray] (3;\k)--(3;\l); 
 \pgfmathtruncatemacro\k{12-\i} 
 \pgfmathtruncatemacro\l{mod(15+\i,18)}   
 \draw[thin,gray] (3;\k)--(3;\l);}    
 
\draw[<->] (-3.5,0) -- (3.5, 0);
%%%%% grey circles at lattice points
\fill [gray] (0;0) circle (2pt);
 \foreach \c in {1,...,3}{%
 \pgfmathtruncatemacro\k{\c*6-1}    
 \foreach \i in {0,...,\k}{% 
   \fill [gray] (\c;\i) circle (2pt);}}  
%%%%% red points in parallelogram
\foreach \x in {0,1, 2, 3, 4, 5}{%
\pgfmathtruncatemacro\y{\x}	
	\filldraw[color=red, fill=red] (1;\y) circle (1mm);}
\fill [red] (2;7) circle (1mm);	
\fill [red] (2;1) circle (1mm);
\fill [red] (0;0) circle (1mm);
\end{tikzpicture}  
\caption{$\la 1, 3, 3 \ra=9$}
\end{subfigure}
\caption{Ternary multiplication includes binary multiplication.}
\label{fig:ternbin}
\end{figure}
\end{center}
\vspace{-3em}
\textbf{Question: Which natural numbers are $3$-prime?}

To answer this question, we need some preliminaries on ternary multiplication. 
\begin{nprop}\label{form1}
The ternary product can be written 
\begin{equation} \label{eq:sym}
\la x, y, z\ra = xy+ yz+ zx -x-y-z +1 .
\end{equation}

\end{nprop}
\begin{proof} We have seen that $\la x, y, 1 \ra =xy$. Increasing the third argument by 1 always adds a strip of $x+y-1$ points along two consecutive edges opposite to those with $x$ and $y$ points (see again Figure~\ref{fig:ternbin}, where the argument increases from 1 to 2).This allows us to conclude the claimed equality,
\begin{equation} \label{eq:form1}
\la x, y, z \ra = xy+ (z-1)(x+y-1) = xy+yz+xz-x-y-z+1.
\end{equation}
\end{proof}
Scholars of symmetric polynomials will recognize (\ref{eq:sym}) as an alternating sum of \textbf{elementary symmetric polynomials},\footnote{This observation can be generalized to construct a family of $n$-ary operations with similar properties.} 
\begin{equation}
\la x,y,z\ra = e_2(x,y,z)-e_1(x,y,z)+e_0(x,y,z).
\end{equation}
However, if you want to mentally compute some ternary products, you may find the formula 
\begin{equation}
 \la x, y, z \ra = xyz -(x-1)(y-1)(z-1) 
 \end{equation}
more convenient.

How can we determine if a number $n$ is $3$-prime? When studying $2$-primes, the first method one usually learns is the \textbf{Sieve of Eratosthenes}, which produces the list of $2$-primes up to a given $N$ by crossing off multiples of those $2$-primes which are at most $\sqrt{N}$. This amounts to eliminating all of the numbers greater than each 2-prime $p$ in the congruence class \hbox{$0 \bmod p$}.  

The proof of Proposition~\ref{form1} indicates how we might sieve for 3-primes. Suppose that $p= x+y-1$ is a $2$-prime, where $x$ and $y$ are natural numbers. We see that 
$$\la x,y,z\ra=xy+(z-1)(x+y-1)=xy+(z-1)p$$ fails to be $3$-prime for all $z\geq 2$, thus we can   also eliminate all of the numbers of the congruence class $xy \bmod p$ which are greater than $xy$ by varying the choice of $z$ in the product $\la x, y, z\ra.$ 

For example, the $2$-prime 3 can be written as 
\[ 3= 3+1-1 \qquad \text{or} \qquad 3=2+2-1. \]
The first case corresponds the choice of $x=3$, $y=1$, which produces the class of ternary products $\la 3, 1, z\ra= 3z$ and eliminates numbers above $3\cdot 1=3$ in the congruence class $0 \bmod 3$ from $3$-primality, as in the usual $2$-primality sieve. But when we take $x=2$, $y=2$, the products of the form 
\[\la 2,2,z\ra= 4+(z-1)(2+2-1)=4+(z-1)3\]
eliminate those numbers that are above 4 and in the congruence class of $4\equiv 1 \bmod 3$ from possible $3$-primality.

We see that for an odd 2-prime $p$, there are $\frac{p+1}{2}$ choices of (unordered) pairs $x$ and $y$ such that \hbox{$p=x+y-1$.} The next proposition shows that each choice produces a distinct congruence class $xy \bmod p$. 

\begin{nprop} \label{prop:distinct} Let $x$ and $w$ be distinct natural numbers between $1$ and a 2-prime $p$, and $w\neq p+1-x$. Then the congruence classes of \hbox{$x(p+1-x)$} and $w(p+1-w)$ modulo $p$ are distinct. 
 \end{nprop}
\begin{proof} We will show the equivalent statement that $x(p+1-x)\equiv w(p+1-w)\bmod p$ implies that $x=w$ or $x=p+1-w$. 
Supposing we have
\[x(p+1-x)\equiv x- x^2 \bmod p\qquad \text{and}\qquad w(p+1-w) \equiv w-w^2 \bmod p,\]
satisfying $x-x^2 \equiv w-w^2 \bmod p$. Then 
\begin{align*}
w^2-x^2-w+x&\equiv 0 \bmod p,\quad \text{ so } \\
(w-x)(w+x-1) &\equiv 0 \bmod p.
\end{align*}
So either $p$ divides $w-x$ or $p$ divides $w+x-1$. Since both $x$ and $w$ are between $1$ and $p$, we have
\[ 1-p\leq w-x \leq p-1 \qquad \text{and} \qquad 1\leq w+x-1\leq 2p-1.\]
Then in the first case, it can only be that $w-x=0$, while in the other case $w+x-1=p$. 
\end{proof}
This has major consequences for how many numbers can be 3-prime! Recall Dirichlet's theorem on arithmetic progressions.
\begin{nthm}[Dirichlet] Let $p$ be a 2-prime and $1\leq k <p$. Then there are infinitely many 2-primes of the form $k+mp$, where $m\in \N$. 
\end{nthm}
A fuller statement of Dirichlet's theorem says that there is the same ``proportion'' of primes in each non-zero congruence class modulo $p$ \cite[Chap. 7]{apostol}. There are $p-1$ such classes for each $p$, and Proposition~\ref{prop:distinct} says that half of them are ruled out from the possibility of 3-primality, in addition to the congruence class $0 \bmod p$. The following lemma further tells us that ruling out congruence classes only needs to happen at the 2-primes -- nothing new is eliminated by ternary products of the form $\la x, y, z \ra$ where $x+y-1$ is 2-composite.

\begin{nlem} \label{lem:reduction} Let $m=ab=x+y-1$. Then there are natural numbers $v$ and $w$ such that $v+w-1=a$ and $xy\equiv vw \bmod a$. Hence if $n=\la x, y, z\ra$, then there exists $z'$ such that $n=\la v,w, z'\ra$.
\end{nlem}
\begin{proof}
We can write 
\[ xy = x(ab+1-x)=xab +x-x^2 \equiv x-x^2 \bmod a.\]
Let $v$ be the least representative of the congruence class $x \bmod a$, and set $w=a+1-v$. Then 
\[vw=v(a+1-v)=va+v-v^2 \equiv v-v^2 \bmod a.\]
Since $v$ and $x$ are in the same congruence class modulo $a$, the claim is proved. 
\end{proof}

To list the 3-primes up to a given $N$, first we can eliminate the $2$-primes below $N$ using the Sieve of Eratosthenes, and now Proposition~\ref{prop:distinct} and Lemma~\ref{lem:reduction} say that we must further eliminate some congruence classes modulo $p$ for some of the 2-primes below $N$. The full method is given in the following ``ternary sieve'' $\mb{TS}$, by proceeding through numerous stages $TS_k$.

\begin{nalgom}[Ternary Sieve]
To determine the $3$-primes less than a given $N$, list the numbers from $2$ to $N$ and perform the following elimination procedure $\mb{TS}$:
\begin{enumerate}
\item $TS_0$: Perform the Sieve of Eratosthenes and create the auxiliary list $\Pi_2(N)$ of $2$-primes at most $N$. 
\item For each $1\leq k \leq \sqrt{\frac{4N-1}{12}}-\frac{1}{2}$ perform elimination step $TS_k$ as follows. Let $T_k=\frac{k(k+1)}{2}$ be the $k$\ts{th} triangular number. For each $p\in \Pi_2(N)$ such that $p\leq \sqrt{N+2T_k}$, eliminate the numbers up to $N$ of the form \hbox{$\la k+1, p-k, p-k\ra+mp$}, for $m\in \N$.
\end{enumerate}
Those numbers that remain among the numbers from 2 to $N$ constitute the list $\Pi_3(N)$ of $3$-primes which are at most $N$.
\end{nalgom}
\begin{proof}
The Sieve of Eratosthenes eliminates the products of the form $\la 1, p ,z\ra$ for $z\geq 2$ by allowing us to add $(z-1)p$ to the product $p\cdot 1=p$. As a small efficiency, starting from the first 2-prime 2 and as $p$ increases towards $\sqrt{N}$, one only needs to cross off products $\la 1,p,z\ra$ with $z\geq p$, as those for $z<p$ will already have been cancelled as products of lesser 2-primes. 

The next stage ($TS_1$) of our sieve for 3-primes eliminates products of the form $\la 2, p-1, z\ra$. Lemma~\ref{lem:reduction} tells us that this can only possibly cross out new numbers if all of the sums
\[ 2+(p-1)-1, \qquad 2+z -1, \qquad\text{and}\qquad (p-1)+z-1 \]
are 2-primes. However, since we will also proceed through this stage using $p$'s in increasing order from among the 2-primes produced by $TS_0$, the elimination of a product will be redundant if $2+z-1$ (which is possibly smallest among the three sums) is a $2$-prime less than $p$. Thus, we can start from $z\geq p-1$ at this stage. Furthermore, this should be done only for those $2$-primes $p$ such that the first possible interesting product
\[\la 2, p-1, p-1 \ra \leq N.\]
Writing 
\[\la 2, p-1, p-1\ra =2(p-1)+(p-2)p=p^2-2,\]
we see that the $TS_1$ stage uses those $2$-primes such that $p\leq\sqrt{N+2}$.

At the $k$\ts{th} stage $TS_k$, we eliminate products of the form $\la k+1, p-k, z\ra$. By considerations similar as in the previous case, this process only needs to happen for $z\geq p-k$, and therefore only for $2$-primes such that (applying formula~\ref{eq:form1})
\[\la k+1, p-k, p-k \ra = (k+1)(p-k)+ (p-k-1)p=p^2-(k^2+k) \leq N.\]
 Rearranging, this condition becomes $p\leq\sqrt{N+2T_k}$.

To obtain the bound on $k$, note that to avoid further unnecessary redundancy we should keep $k+1\leq p-k$ in the product $\la k+1, p-k, z\ra$.\footnote{Otherwise, once $k+1$ becomes bigger than $p-k$ we start transiting the same choices of pairs for the first two inputs, but in the opposite direction.} Since $z\geq p-k$ during $TS_k$, the largest possible value of $k$ must satisfy
\[ \la k+1,k+1,k+1\ra = 3k^2+3k+1\leq N.\]
Completing the square and solving the inequality, one finds that 
\begin{equation}
k\leq \sqrt{\frac{4N-1}{12}}-\frac{1}{2}
\end{equation}
is a sufficient bound. Note that by Lemma~\ref{lem:reduction}  $\la k+1, k+1, k+1\ra$ will only eliminate a new congruence class if $p=2k+1$ is 2-prime.
\end{proof}

This algorithm is not hard to implement on a computer, and a search for the 3-primes up to 10,000,000 reveals a very short list. 
\[ 2,\ 3,\ 5,\ 11,\ 17,\ 41 \]

At OEIS (\href{https://oeis.org/A014556}{A014556}) one learns that these are ``Euler's `lucky' numbers,'' those 2-primes $p$ such that 
\begin{equation}
n^2-n+p
\end{equation}
is 2-prime for all $1\leq n \leq p-1$. This confirms that $3$-primes are somehow ``extra'' prime, but these numbers are significant for a deeper reason. We might add the number $1$ to Euler's list, as it vacuously satisfies the defining condition, so obtaining an ``augmented lucky numbers'' list. The augmented list is then exactly the set of integers $k$ such that $4k-1$ is a \textbf{Heegner number}. The full list of Heegner numbers is 
\[1,2,3,7,11,19,43,67,163.\]
We next explain their significance.

\section{The Class Number One Problem}
Recall that a \textbf{quadratic number field} $\Q(\sqrt{n})$ is an extension of the rational numbers $\Q$ obtained by adjoining the root of an irreducible degree-two polynomial. Just as the integers $\Z$ sit inside the rationals, each quadratic number field has its own set of integers.

\begin{ndefn} The \textbf{ring of integers} of a quadratic number field $K=\Q(\sqrt{n})$ is the subset of elements which are roots of some monic polynomial with coefficients in $\Z$. It is denoted $\mc{O}_K$.
\end{ndefn}

Classic examples include the \emph{Gaussian integers } $\Z[i]$ inside $\Q(i)$, and the \emph{Eisenstein integers} $\Z[\omega]$ inside $\Q(\sqrt{-3})$, where $\omega=-\frac{1}{2}+ \frac{\sqrt{-3}}{2}$. In general one has the following uniform description of rings of \emph{quadratic integers} \cite[189]{irelandRosen}.
\begin{equation}\label{eq:integers}
\mc{O}_K= \begin{cases} \Z[\sqrt{n}] & \text{if } n \not\equiv 1 \bmod 4\\
\Z\lt[\frac{-1+\sqrt{n}}{2}\rt] & \text{if } n \equiv 1 \bmod 4 .
\end{cases} 
\end{equation}

It is well-known that both the Gaussian integers and the Eisenstein integers admit unique factorization into irreducible elements, just as the ordinary integers admit unique factorization into 2-primes. But other rings of quadratic integers do not. To cite a common example, in $\Z[\sqrt{-5}]$ the number $6$ admits decompositions as both $2\cdot 3$ and $(1+\sqrt{-5})(1-\sqrt{-5})$. This leads to the following natural question.

\begin{center}
For which $n$ does the ring of integers of $\Q(\sqrt{n})$ have unique factorization?
\end{center}
 
Beyond the description of (\ref{eq:integers}), there is a bifurcation in the approach to this question according to whether $n$ is positive or negative; that is, whether the quadratic field is \emph{real} or \emph{imaginary}. These two types are of extremely different character. For instance, there are very few units (invertible elements) in the ring of integers of imaginary fields, while there are infinitely many in the real case. We are concerned here with the imaginary case, for which a complete answer to the question above is known. 

\begin{nthm} \label{thm:cno1} For a natural number $n$, the ring of integers of $\Q(\sqrt{-n})$ has unique factorization if and only if $n$ is a Heegner number: 1, 2, 3, 7, 11, 19, 43, 67, or 163.
\end{nthm}

This theorem, has a long, interesting history with origins in the study of quadratic forms going back to Fermat, Lagrange, Legendre, Gauss, etc.\footnote{See, for instance, \cite{goldfeldGauss}.} The answer was guessed by Gauss and was proved by Heegner in 1952, but this proof was only accepted by the mathematical community after Heegner's death and the appearance of proofs in the 1960's by established mathematicians Baker and Stark. Moreover, the answer to the unique factorization problem is just one part of a larger problem called \textbf{Gauss' class number problem}, resolved by Goldfeld-Gross-Zagier in 1985. Theorem~\ref{thm:cno1} above addresses just the \emph{class number one} problem, with classes referring to equivalence classes either of ideals in $\mc{O}_K$ or of a related set of quadratic forms.  When the class number of $\mc{O}_K$ is one, it implies that $\mc{O}_K$ has the unique factorization property; for background, see \cite{coxPrimes}.

Within the class number one problem, the two cases of (\ref{eq:integers}) are treated differently. Remember that a quadratic number field is obtained by adjoining to $\Q$ a root $\alpha$ of some polynomial \hbox{$ax^2+bx+c$}, which root has formula
\[\alpha=\frac{-b\pm\sqrt{b^2-4ac}}{2a}.\]
Since we are adjoining $\alpha$ to the rational numbers, everything other than $\sqrt{b^2-4ac}$ can be disregarded, and in fact the discriminant $D=b^2-4ac$ determines the number field. Since $D\equiv b^2$ modulo 4 it can only be congruent to $0$ or $1$. When $D\equiv 0 \bmod 4$, there is a factor of 4 that can be pulled out of the square root so that 
\[\Q(\sqrt{D})=\Q\lt(\sqrt{\frac{D}{4}}\rt).\]
This is to say that for a quadratic field $\Q(\sqrt{n})$, $n$ is usually understood to be square-free, although it may come by adjoining the root of a polynomial with discriminant $D=4n$. This motivates the following definition.

\begin{ndefn} The \textbf{field discriminant} $d_K$ of $K=\Q(\sqrt{n})$ is  
\[ d_K= \begin{cases}
n & \text{if } n \equiv 1 \bmod 4 \\
4n & \text{otherwise.}
\end{cases}
\]
\end{ndefn}

Returning to the list of Heegner numbers, we see that $\Q(\sqrt{-1})$ and $\Q(\sqrt{-2})$ are the only cases where $d_K\equiv 0 \bmod 4$. That is, $-n\equiv 1 \bmod 4$ for every other Heegner number $n$. Often, Theorem~\ref{thm:cno1} is stated by giving instead the list of negative field discriminants $D$ such that $\Q(\sqrt{D})$ has \textbf{class number} $h(D)$ equal to one. Then, instead of the Heegner numbers, we have the slightly modified list
\[D= -3, -4, -7, -8, -11, -19, -43, -67, -163.\]

In 1902, Landau was able to prove that $\Q(\sqrt{-1})$ and $\Q(\sqrt{-2})$ are the only imaginary quadratic fields with even (divisible by 4, really) discriminant and unique factorization.\footnote{Actually, he proved a slightly broader statement in terms of quadratic forms; see \cite[{Theorem 2.18}]{coxPrimes}. } The proof of this fact is quite elementary, but the proofs of Heegner, Baker and Stark that cover the odd discriminant case require much more sophistication.\footnote{Heegner's and Stark's proofs use modular forms, while Baker's approach involves bounds on logarithms of linear forms of algebraic numbers.} In 1913, Rabinowitsch provided another elementary characterization of the odd case.

\begin{nthm}[{\cite{rabinowitsch}}]\label{thm:rab} Let $D<0$ and $D\equiv 1 \bmod 4$. Then 
\[x^2-x+ \frac{1+\ab{D}}{4}\quad\text {is prime for each}\quad x= 1,2,\dots, \frac{\ab{D}-3}{4}\]
if and only if the integers of the field $\Q(\sqrt{D})$ admit unique factorization.
\end{nthm}

This theorem does not appear to have been directly useful for solving the class number one problem, but it does link it to the list of Euler's lucky primes, and so to the list of 3-primes. Just as we augmented Euler's lucky primes by adding 1, from now on we will consider 1 as a 3-prime in the sense that it is not representable by a non-degenerate hexagonal or parallelogrammatic configuration. We then see that the augmented lucky numbers/known 3-primes
\[1,2,3,5,11,17, 41\] 
account for all of the negative odd discriminants of class number one,
\[-3, -7, -11, -19, -43, -67, -163.\]
Next we show that the 3-primes known from the ternary sieve exactly coincide with the augmented lucky numbers.
\begin{nthm} \label{thm:3prime} A number is $3$-prime if and only if it is among the augmented lucky numbers.
\end{nthm}
\begin{proof}
If $n$ is not 3-prime, then $n$ dots can be arranged into a lattice parallelogram or hexagon such that two distinct pairs of sides have at least two points along the edge. When this is the case, either $n$ is already $2$-composite, in which case a parallelogrammatic representation exists, or not, in which case a true hexagonal representation exists. Supposing the latter is the case, the hexagon can be ``completed'' to a parallelogram by adding two triangles along opposite edges. If the the sides abutting these triangles contain $k$ points, then the completed parallelogram will have $n+2T_{k-1}$ points (see Figure~\ref{fig:completion}). 

On the other hand, if a number $n$ is $3$-prime, then the only representation it admits is a row of $n$ dots. In other words, the smallest triangles that can be adjoined in order to obtain a parallelogram are those of size $T_{n-1}$ (see Figure~\ref{fig:5prime}). Equivalently, $n+2T_k$ is $2$-prime for $k=1,2,\dots,{n-1}$, as is $n$ itself. 

\begin{center}
\begin{figure}[htbp!]
\begin{subfigure}[b]{0.45\textwidth}
\begin{tikzpicture}[scale=0.8]
%%%  define vertices with coordinates
\coordinate (0;0) at (0,0); 
\foreach \c in {1,...,4}{%  
\foreach \i in {0,...,5}{% 
\pgfmathtruncatemacro\j{\c*\i}
\coordinate (\c;\j) at (60*\i:\c);  
} }
\foreach \i in {0,2,...,10}{% 
\pgfmathtruncatemacro\j{mod(\i+2,12)} 
\pgfmathtruncatemacro\k{\i+1}
\coordinate (2;\k) at ($(2;\i)!.5!(2;\j)$) ;}

\foreach \i in {0,3,...,15}{% 
\pgfmathtruncatemacro\j{mod(\i+3,18)} 
\pgfmathtruncatemacro\k{\i+1} 
\pgfmathtruncatemacro\l{\i+2}
\coordinate (3;\k) at ($(3;\i)!1/3!(3;\j)$)  ;
\coordinate (3;\l) at ($(3;\i)!2/3!(3;\j)$)  ;
 }
 
 \foreach \i in {0,4,...,20}{% 
\pgfmathtruncatemacro\j{mod(\i+4,24)} 
\pgfmathtruncatemacro\k{\i+1} 
\pgfmathtruncatemacro\l{\i+2}
\pgfmathtruncatemacro\m{\i+3}
\coordinate (4;\k) at ($(4;\i)!1/4!(4;\j)$)  ;
\coordinate (4;\l) at ($(4;\i)!2/4!(4;\j)$)  ;
\coordinate (4;\m) at ($(4;\i)!3/4!(4;\j)$)  ;
 }

%%shape
\filldraw[green!20!white] (2;8) -- (2;10) -- (2;0) -- (2;2) -- (2;4)-- (2;6) ;
\filldraw[yellow!20!white] (2;0) -- (2;2) -- (4;2)
				 (2;6) -- (2;8)-- (4;14) ;
 %%%%%%%%% draw lines %%%%%%%%
 % positive slope
 \foreach \i in {0,...,8}{% 
 \pgfmathtruncatemacro\k{\i}
 \pgfmathtruncatemacro\l{20-\i}
 \draw[thin,gray] (4;\k)--(4;\l);
 % negative slope
 \pgfmathtruncatemacro\k{12-\i} 
 \pgfmathtruncatemacro\l{mod(16+\i,24)}   
 \draw[thin,gray] (4;\k)--(4;\l); 
 % horizontal
 \pgfmathtruncatemacro\k{16-\i} 
 \pgfmathtruncatemacro\l{mod(20+\i,24)}   
 \draw[thin,gray] (4;\k)--(4;\l);}    
 %axes
 \draw[<->] (-4.5,0) -- (4.5, 0);
 \draw[<->] (0,-4.5) -- (0, 4.5);
%%% draw gray circles at lattice points
\fill [gray] (0;0) circle (2pt);
 \foreach \c in {1,...,4}{%
 \pgfmathtruncatemacro\k{\c*6-1}    
 \foreach \i in {0,...,\k}{% 
   \fill [gray] (\c;\i) circle (2pt);}}  
%%%%% points in hexagon
\foreach \x in {0,1, 2, 3, 4, 5}{%
\pgfmathtruncatemacro\y{\x}	
	\filldraw[color=red, fill=red] (1;\y) circle (1mm);}
\foreach \x in {0,...,11}	{ \pgfmathtruncatemacro\y{\x}
\fill [red] (2;\x) circle (1mm);}	
\fill [red] (0;0) circle (1mm);
%triangles
\fill [purple] (3;1) circle (1mm)
	[purple] (3;2) circle (1mm)
	[purple] (4;2) circle (1mm);
\fill [purple] (3;11) circle (1mm)
	[purple] (3;10) circle (1mm)
	[purple] (4;14) circle (1mm);	
\end{tikzpicture}  
\caption{$\la3,3,3\ra=19$ is not 3-prime because \hbox{$19+2T_2=19+6=25$} is 2-composite. }
\label{fig:completion}
\end{subfigure}
\begin{subfigure}[b]{0.45\textwidth}
\begin{tikzpicture}[scale=0.8]
%%%  define vertices with coordinates
\coordinate (0;0) at (0,0); 
\foreach \c in {1,...,4}{%  
\foreach \i in {0,...,5}{% 
\pgfmathtruncatemacro\j{\c*\i}
\coordinate (\c;\j) at (60*\i:\c);  
} }
\foreach \i in {0,2,...,10}{% 
\pgfmathtruncatemacro\j{mod(\i+2,12)} 
\pgfmathtruncatemacro\k{\i+1}
\coordinate (2;\k) at ($(2;\i)!.5!(2;\j)$) ;}

\foreach \i in {0,3,...,15}{% 
\pgfmathtruncatemacro\j{mod(\i+3,18)} 
\pgfmathtruncatemacro\k{\i+1} 
\pgfmathtruncatemacro\l{\i+2}
\coordinate (3;\k) at ($(3;\i)!1/3!(3;\j)$)  ;
\coordinate (3;\l) at ($(3;\i)!2/3!(3;\j)$)  ;
 }
 
 \foreach \i in {0,4,...,20}{% 
\pgfmathtruncatemacro\j{mod(\i+4,24)} 
\pgfmathtruncatemacro\k{\i+1} 
\pgfmathtruncatemacro\l{\i+2}
\pgfmathtruncatemacro\m{\i+3}
\coordinate (4;\k) at ($(4;\i)!1/4!(4;\j)$)  ;
\coordinate (4;\l) at ($(4;\i)!2/4!(4;\j)$)  ;
\coordinate (4;\m) at ($(4;\i)!3/4!(4;\j)$)  ;
 }

%%shape
\filldraw[yellow!10!white] (2;0) -- (4;6) -- (2;6) -- (4;18);
 %%%%%%%%% draw lines %%%%%%%%
 % positive slope
 \foreach \i in {0,...,8}{% 
 \pgfmathtruncatemacro\k{\i}
 \pgfmathtruncatemacro\l{20-\i}
 \draw[thin,gray] (4;\k)--(4;\l);
 % negative slope
 \pgfmathtruncatemacro\k{12-\i} 
 \pgfmathtruncatemacro\l{mod(16+\i,24)}   
 \draw[thin,gray] (4;\k)--(4;\l); 
 % horizontal
 \pgfmathtruncatemacro\k{16-\i} 
 \pgfmathtruncatemacro\l{mod(20+\i,24)}   
 \draw[thin,gray] (4;\k)--(4;\l);}    
 %axes
 \draw[<->] (-4.5,0) -- (4.5, 0);
 \draw[<->] (0,-4.5) -- (0, 4.5);
%%% draw gray circles at lattice points
\fill [gray] (0;0) circle (2pt);
 \foreach \c in {1,...,4}{%
 \pgfmathtruncatemacro\k{\c*6-1}    
 \foreach \i in {0,...,\k}{% 
   \fill [gray] (\c;\i) circle (2pt);}}  
%%%%% points in hexagon
\fill [red] (2;0) circle (1mm)
		(2;6) circle (1mm)
		(1;0) circle (1mm)
		(1;3) circle (1mm);	
\fill [red] (0;0) circle (1mm);
%triangles
\fill [purple] (2;1) circle (1mm) (2;2) circle (1mm) (2;3) circle (1mm) (2;4) circle (1mm) (2;5) circle (1mm) (2;7) circle (1mm) (2;8) circle (1mm) (2;9) circle (1mm)(2;10) circle (1mm)(2;11) circle (1mm);
\fill [purple] (1;1) circle (1mm) (1;2) circle (1mm) (1;4) circle (1mm) (1;5) circle (1mm);
\fill [purple] (3;4) circle (1mm) (3;5) circle (1mm) (3;13) circle (1mm) (3;14) circle (1mm);
\fill [purple] (4;18) circle (1mm) (4;6) circle (1mm);
\end{tikzpicture}  
\caption{5 is 3-prime because two $T_4$ triangles are the smallest that can be added to reach a $2$-composite. }
\label{fig:5prime}
\end{subfigure}
\caption{Relating $3$-factorizations and $2$-factorizations.}
\end{figure}
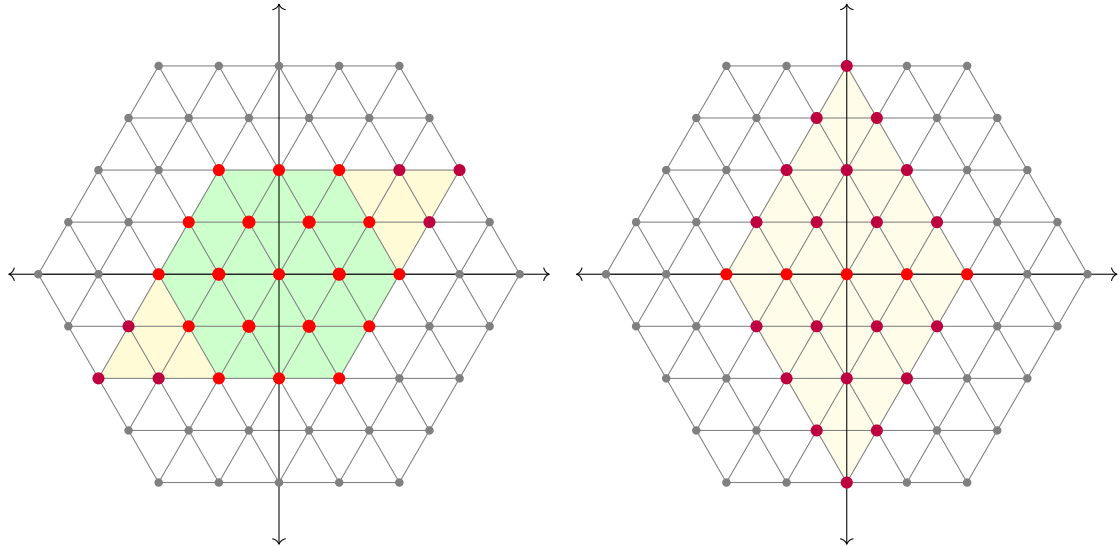
\end{center}
\vspace{-3em}

Now examine the polynomial $x^2-x+n$, and observe that $x^2-x=x(x-1)=2T_{x-1}$ when $x$ is a natural number. Then, the condition that $x^2-x+n$ is 2-prime for all $x$ from 1 to $n-1$ is equivalent to statement that every number in the set
\[ \{n, n+2T_1, n+2T_2, \dots n+2T_{n-1} \} \]
is $2$-prime. This is plainly equivalent to the characterization of $3$-primality just given.
\end{proof}
Invoking the theorems of Rabinowitsch and Heegner/Baker/Stark, we conclude the following.
\begin{ncor} There are only finitely many 3-primes. Including 1, these are 1, 2, 3, 5, 11, 17, and 41. \label{cor:3primes}
\end{ncor}
\begin{proof} We see that if there were $3$-primes beyond the list of augmented lucky numbers, they would give imaginary quadratic fields with unique factorization and odd discriminant. But there are only seven of these from the solution of the class number one problem.
\end{proof}

This is a very heavy-handed proof, especially compared to common proofs of the finitude of 2-primes. While it would be extremely desirable to find a proof that relied only upon ternary multiplication, the historical difficulty of the class number one problem suggests that this might be out of reach.

\section{Applications}
Besides determining which numbers are 2-prime and which are 2-composite, one of the most basic questions one can ask in number theory is how to determine the factorization of numbers which are 2-composite. The proof of Theorem~\ref{thm:3prime} can be retooled to produce \textbf{$3$-factorizations} of natural numbers, by which we mean representations of a number as a ternary product.

For example, by Corollary~\ref{cor:3primes}, we know that 19 is 3-composite. To find its 3-factorizations, we can add twice a triangular number to 19 to see if we obtain a $2$-composite number. Then, since we know that a $2$-composite of the form $19+2T_k$ can be represented by a parallelogram, there is an equation $19+2T_k=ab$ where neither of $a$ and $b$ is equal to 1. Removing the corner triangles (consisting of $2T_k$ points) from this parallelogram, we get a hexagon whose sides give a non-trivial 3-factorization of 19. 

Note that there may be several $k$ for which $19+2T_k$ is 2-composite and multiple parallelograms that represent each of those 2-composites. For instance, $19+2=21=7\cdot 3$. Removing two points from the opposite corners of a $7$ by $3$ parallelogram, we get a hexagon with pairs of sides of lengths 2, 2 and 6, so $19= \la 2, 2, 6\ra$. But $19+2T_2=19+6=25= 5\cdot 5$ as well, and removing the $T_2$ triangles from the $5$ by $5$ parallelogram gives us the 3-factorization $19=\la 3,3,3 \ra$ of Figure~\ref{fig:completion}. In general, we have the following.

\begin{nprop} If $n+2T_k=ab$ for $a,b > k$, then $n=\la a-k, b-k,k+1 \ra$.
\end{nprop}
\begin{proof}
Construct a lattice parallelogram which has $a$ and $b$ points along opposite pairs of edges. Since $a$ and $b$ are bigger than $k$, we can remove lattice triangles with $k$ points along each edge from opposite corners of the parallelogram. If $a$ and $b$ are $k+1$, then removal of these triangles yields $n$ points in a row and the factorization $n=\la 1, 1, n \ra $. In case exactly one of $a$ or $b$ is $k+1$ (assume it is $a$), the removal produces a new parallelogram and the factorization $n=\la 1, b-k, k+1\ra$. Otherwise this removal produces a true lattice hexagon. The number of points in opposite pairs of edges of this hexagon are $a-k$, $b-k$, and $k+1$ which yields the factorization $n=\la a-k, b-k, k+1\ra$.
\end{proof}

We see that the proposition also covers ``degenerate'' 3-factorizations which are either trivial ($n=\la1,1,n\ra$) or reduce to binary products, though we are most interested in the 3-factorizations where each factor is at least 2. One could obtain all of these hexagonal representations of $n$ as follows.

 Suppose $n$ has 3-factorization \hbox{$n=\la x,y,z \ra$}, where $z$ is the least among the three factors. This 3-factorization can be discovered by examining parallelograms which represent $n+2T_{z-1}$ and removing the triangles in opposite corners.  The smallest $3$-factor of $n$ is as large as possible when $n=\la z, z, z\ra$, meaning that to recover all 3-factorizations, one needs to examine the $2$-factorizations of all the numbers $n+2T_k$ for $1\leq k \leq z-1$, where $z$ is the largest number satisfying 
\[\la z,z,z\ra =3z^2-3z+1\leq n.\]

Since ternary multiplication results in a number system with finitely many $3$-primes, the fact that many numbers admit multiple 3-factorizations is not surprising. The question of exactly how many distinct 3-factorizations (up to reordering the factors) a number admits, and how this statistic may be further related to the class numbers of quadratic fields could be of interest for future research. 

 \begin{table}[htbp]
\caption{Number of 3-factorizations of small natural numbers}
\centering 
\begin{tabular}{c|c|c|c|c|c|c|c|c|c|c|c|c|c|c|c|c|c|c|c|c}
\toprule
$n$ &1 & 2 & 3 & 4 & 5 & 6 & 7 & 8 & 9 & 10 & 11 & 12 & 13 & 14 & 15 & 16 & 17 & 18 & 19 & 20 \\\hline
3-factorizations &1 & 1 & 1 & 2 & 1 & 2 & 2 & 2 & 2 & 3 & 1 & 3 & 2 & 3 & 2 & 4 & 1 & 4 & 3 & 3\\\bottomrule
\end{tabular}
\end{table}

We will close this discussion with a few applications of this line of thinking to elementary number theory. The first is a 2-primality test that comes from the following partial converse of Proposition~\ref{prop:distinct}. 

\begin{nprop} \label{prop:conv}Let $n=pr$ be an odd $2$-composite number and $p,r\geq 3$.  Then there are distinct $x$ and $w$ where $1\leq x,w \leq \frac{n+1}{2}$ such that $x(n+1-x)\equiv w(n+1-w) \bmod n$. 
\end{nprop}
\begin{proof} 
We can assume $p\leq r$ by choosing $p$ to be the smallest 2-prime factor of $n$.  We will show that the claim is true for some $x\leq \frac{n+1}{2}$ and $w=x+p\leq\frac{n+1}{2}$, though the statement may be true for other choices of $x$ and $w$ as well. In order to obtain 
\[ x(n+1-x) \equiv w(n+1-w) \bmod n\]
we need
\begin{align*}
x(n+1-x) &\equiv (x+p)(n+1-x-p) \bmod n, \quad \text{so} \\
nx+x-x^2 & \equiv nx+x-x^2-px+np+p-px-p^2 \bmod n. \\
\end{align*}
Subtracting and collecting terms, we have
\[2px+p^2-p=p(2x+p-1)\equiv 0 \bmod n.\]
This is satisfied if and only if $2x+p-1\equiv 0 \bmod r$. By varying $x$, we can arrange $2x+p-1$ to take the value of any even number from 
$$p+1 \qquad \text{to}\qquad  2\lt(\frac{n+1}{2}-p\rt)+p-1=pr+1-2p+p-1=pr-p.$$
By showing that the even number $2r$ lies in this range, we will establish the existence of $x$ and $w$.
First, $p+1\leq 2r$ because $r\geq p$ and both are at least 3. Next, the inequality $2r\leq pr-p$ holds if and only if 
\begin{align*}
pr-2r-p &\geq 0\\
(p-2)r-p\geq 0
\end{align*}
which holds because $p\geq 3$ and $r\geq p$. 
\end{proof}

The statement of this proposition seems to be true for any $r\geq 2$, and so for every $2$-composite that is not a pure power of $2$ rather than just for odd $n$. However, 2-primality tests usually are usually only meant for odd numbers since even numbers can be tested instantly, so we ignore this other case. Thus we have the following 2-primality test, which is an immediate consequence of Propositions~\ref{prop:distinct} and~\ref{prop:conv}.

\begin{nthm} \label{thm:2test} Let $n$ be an odd natural number. Then $n$ is 2-prime if and only if the congruence classes of $x(n+1-x) \bmod n$ are distinct for every $x$ between $1$ and $\frac{n+1}{2}$.
\end{nthm}

This test doesn't appear to be very efficient -- as stated it requires about half as many computations as the size of the number $n$.\footnote{The best (deterministic) $2$-primality testing algorithms require a number of computations which is polynomial in $\log n$, instead of linear in $n$ as this one is. } However a simple observation makes it slightly more suitable for hand calculation with small numbers. 
\begin{nlem} Let $T_k=\frac{k(k+1)}{2}$ for any $k=0, 1, 2, 3,\dots$ and let $n$ and $x$ be natural numbers with $1\leq x \leq n$. Then \hbox{$x(n+1-x) \equiv -2T_{x-1} \bmod n$.}
\end{nlem}
\begin{proof} Recall that $x(n+1-x)\equiv x-x^2 \bmod n$ and $x-x^2=x(1-x)=-2T_{x-1}$.
\end{proof}

This lemma makes it quite easy to write down the congruence classes of interest for Theorem~\ref{thm:2test}. To make the list, start from $x=1$, in which case $x(n+1-x)\equiv-2T_0\equiv0\bmod n$, and then to get from $-2T_x$ to $-2T_{x+1}$, just subtract $2(x+1)$. We illustrate this now for $n=15$:
\[ 15\equiv \tcp{0}\bmod 15 \xra{\text{subtract 2}} 13 \xra{\text{subtract 4}} 9 \xra{\text{subtract 6}} \tcb{3} \xra{\text{subtract 8}} -5\equiv 10 \xra{\text{subtract 10}} \tcp{0} \xra{\text{subtract 12}} \tcb{3} \xra{\text{subtract 14}} {4} .\]

We see that the appearance of the congruence class $3 \bmod 15$ at $x=4$ and $x=7$ indicates that 15 is $2$-composite. Also notice the coincidence at $x=1$ and $x=6$, indicating that the type of repetition produced in the proof of Proposition~\ref{prop:conv} occurs not only at intervals of length equal to the smallest prime $p$. 

The coincidence of congruence classes $-2T_k\equiv -2T_l \bmod n$ means that $2T_l-2T_k$ is a multiple of $n$.  The next proposition goes further, relating the value $l-k$ to the $2$-factorization of $n$.

\begin{nprop}
With notation as before, if $2(T_l-T_k)=mn$ for $0\leq k,l \leq \frac{n-1}{2}$, distinct and $n>3$, then both of the pairs $(l-k, n)$ nor $(l+k+1, n)$ have greatest common divisor (gcd) greater than 1.
\end{nprop}
\begin{proof}
We have 
\[2({T_l}-T_k)=l^2+l-k^2-k=(l-k)(l+k+1)=mn.\]
If $\gcd(l-k,n)=1$, then $l-k$ divides $m$ and $l+k+1=sn$ for some integer $s$. But since $0\leq k,l \leq \frac{n-1}{2}$, and $k$ and $l$ are distinct, we have
\[2\leq l+k+1\leq n-1,\]
so $l+k+1=sn$ is impossible.

On the other hand, if $\gcd(l+k+1,n)=1$, then $l-k=tn$ for some integer $t$. But $-\frac{n}{2}<l-k<\frac{n}{2}$, so the only possibility is $t=0$ contradicting that $l$ and $k$ are distinct. 
\end{proof}

Taking the gcd of natural numbers can be done efficiently by Euclid's algorithm. Thus, the major cost in the following 2-factorization algorithm is generating the list of congruence classes $x(n+1-x) \equiv -2T_{x-1} \bmod n$.

\begin{nthm}[2-factorization algorithm]\label{thm:2fact} Given an odd natural number $n$, one can obtain a non-trivial $2$-divisor of $n$ as follows.
\begin{enumerate}
\item List the congruence classes $-2T_k \bmod n$ for $k=0,1,2 \dots$ until there is a repetition $$-2T_k\equiv -2T_l \bmod n, \quad k\neq l$$ In case there is no repetition up through $k=\frac{n-1}{2}$, then conclude $n$ is 2-prime. 
\item Otherwise compute either $\gcd(l-k, n)$ or $\gcd(l+k+1,n)$. The output will be a non-trivial divisor $d$ of $n$. 
\end{enumerate}
Steps 1 and 2 can then be iterated on $d$ and $n_1=\frac{n}{d}$ to obtain a complete $2$-factorization of $n$.
\end{nthm}

One may recognize in this algorithm a formal similarity with \emph{Pollard's rho algorithm}, which also finds a non-trivial factor of $n$ by taking the gcd of numbers after finding a repetition in a sequence. However the discovery of a repetition in the algorithm of Theorem~\ref{thm:2fact} does not mean that there is a ``cycle'' in the sequence as it does in Pollard's algorithm.

 Closer examination reveals that this algorithm actually has more in common with the \emph{Fermat factorization method} which finds factors of $n$ by representing it as a difference of squares, $n=a^2-b^2$.\footnote{See \cite[Ch. 5]{bressoud} or \cite[Ch. 3]{hpsCrypto}, for instance, for descriptions of these other algorithms.} To see this, let $u=\frac{n+1}{2}$, so that when $x=u$ the product $x(n+1-x)$ is $u^2$. Then all of the other products in the list are \hbox{$(u+a)(u-a)=u^2-a^2$} for some $a$. In seeking a match
\[u^2-a^2 \equiv u^2-b^2 \bmod n,\]
we are really seeking a solution to $a^2-b^2\equiv 0 \bmod n$, or $a^2-b^2=mn$ for some integer $m$. 

The ideas of the Fermat factorization method form the basis of the fastest known integer factorization algorithms, the \emph{quadratic sieve} and \emph{general number field sieve} \cite{pomerance}.  It remains to be seen whether the algorithm of Theorem~\ref{thm:2fact} admits improvements that could make it competitive. For now, it is a curiousity which we hope encourages the reader to explore the plunder of ideas which may come from non-binary thinking.

\textbf{Acknowledgements:} I thank an anonymous referee for helpful comments and criticisms. I'm also grateful to Craig Jensen for permission to think about this topic as a master's student, Ken Holladay for a helpful suggestion, Padi Fuster and Bradford Fournier for encouragement, and Mahir Can and Joseph Silverman for comments on drafts of this manuscript.

\printbibliography

\end{document}